\documentclass[11pt]{article}
\usepackage{amssymb,amsfonts,amsmath,latexsym,epsf,tikz,url}
\usepackage{graphicx}
\usepackage[font={small}]{caption}
\usepackage{cite}

\newtheorem{theorem}{Theorem}[section]

\newtheorem{lemma}[theorem]{Lemma}

\newtheorem{example}[theorem]{Example}
\newcommand{\proof}{\noindent{\bf Proof.\ }}
\newcommand{\qed}{\hfill $\square$\medskip}

\begin{document}

\title{\textbf{Extremal trees for Maximum Sombor index with given degree sequence}}
\author{F. Movahedi$^{a}$\footnote{Corresponding author \, E-mail: f.movahedi@gu.ac.ir}}

\maketitle

\begin{center}
$^a$ Department of Mathematics, Faculty of Sciences, Golestan University, Gorgan, Iran.
\end{center}
\maketitle

\begin{abstract}
Let $G=(V, E)$ be a simple graph with vertex set $V$ and edge set $E$. The Sombor index of the graph $G$ is a degree-based topological index, defined as
$$SO(G)=\sum_{uv \in E}\sqrt{d(u)^2+d(v)^2},$$
in which $d(x)$ is the degree of the vertex $x \in V$ for $x=u, v$. \\
In this paper, we characterize the extremal trees with a given degree sequence that maximizes the Sombor index. 
 	
\end{abstract}

\medskip
\noindent{\bf AMS Subj.\ Class.:} 05C35, 05C90.

\section{Introduction}

\indent In \cite{1}, Gutman defined a new vertex degree-based topological index, named the Sombor index, and defined for a graph $G$ as follows

$$SO(G)=\sum_{uv \in E(G)}\sqrt{d(u)^2+d(v)^2},$$

\noindent where $d(u)$ and $d(v)$ denote the degree of vertices $u$ and $v$ in $G$, respectively. \\

Other versions of the Sombor index are induced and studied in \cite{1, 2, 3, 4, 5}. Guman \cite{1} showed that the Sombor index is minimized by the path and maximized by the star among general trees of the same size. In \cite{a4} the extremal values of the Sombor index of trees and unicyclic graphs with a given maximum degree are obtained. Deng et al. \cite{Deng} obtained a sharp upper bound for the Sombor index and the reduced Sombor index among all molecular trees with fixed numbers of vertices, and characterized those molecular trees achieving the extremal value. In \cite{Li} characterized the extremal graphs with respect to the Sombor index among all the trees of the same order with a given diameter. R\'eti et al. \cite{Ret} characterized graphs with the maximum Sombor index in the classes of all connected unicyclic, bicyclic, tricyclic, tetracyclic, and pentacyclic graphs of a given order. In this paper, we focus on the following natural extremal problem of Sombor index. \\

\vspace{0.5cm}

\noindent {\bf Problem 1.} Find extremal trees of Sombor indices with a given degree sequence and characterize all extremal trees which attain the extremal values.

\vspace{0.5cm}

Let $T=(V, E)$ be a simple and undirected tree with vertex set $V(G)=\{ v_1, \ldots, v_n\}$ and the edge set $E(G)=\{e_1, \ldots, e_m\}$. The set $N_T(u)=\{v \in V \vert uv \in E\}$ is called the neighborhood of vertex $u \in V$ in tree $T$. The number of edges incident to vertex $u$ in $G$ is denoted $d(u)=d_u$. A leaf is a vertex with degree $1$ in tree $T$. The minimum degree and the maximum degree of $T$ are denoted by $\delta$ and $\Delta$, respectively. The distance between vertices $u$ and $v$ is the minimum number of edges between $u$ and $v$ and is denoted by $d(u, v)$. The degree sequence of the tree is the sequence of the degrees of non-leaf vertices arranged in non-increasing order. Therefore, we consider $(d_1, d_2, \ldots, d_k)$ as a degree sequence of the tree $T$ where $d_1\geq d_2\geq \cdots \geq d_k \geq 2$. A tree is called a maximum optimal tree if it maximizes the Sombor index among all trees with a given degree sequence. \\

In this paper, we investigate the extremal trees which attain the maximum Sombor index among all trees with given degree sequences.

\section{Preliminaries}
In this section, We prove Some lemmas that are used in the next main results. \\

\begin{lemma}\label{lemma1}
For function $g(x, y)=\sqrt{x^2+y^2}$, if $x\leq y$ then $g(x, 1)\leq g(y, 1)$.
\end{lemma}
\proof
If $x\leq y$, then $x^2+1\leq y^2+1$ and consequently, $\sqrt{x^2+1}\leq \sqrt{y^2+1}$. Therefore, $g(x, 1)\leq g(y, 1)$.
\qed

\begin{lemma}\label{lemma2}
Let $f(x)=\sqrt{x^2+a^2}-\sqrt{y^2+b^2}$ with $a, b, x\geq 1$. Then $f(x)$ is an increasing function for every $a\leq b$ and a decreasing function for every $a>b$. 
\end{lemma}
\proof
We have that
$$f'(x)=\frac{x}{\sqrt{x^2+a^2}}-\frac{x}{\sqrt{x^2+b^2}}.$$
We consider function $\hat{f}(y)=\frac{x}{\sqrt{x^2+y^2}}$ where $y\geq 1$. The derivative of function $\hat{f}(y)$ is $\hat{f}'(y)=\frac{-xy}{(x^2+y^2)\sqrt{x^2+y^2}}<0$. Therefore, $\hat{f}(y)$ is a decreasing function for every $y\geq 1$. Hence, if $a\leq b$, $\frac{x}{\sqrt{x^2+a^2}}=\hat{f}(a)\geq \hat{f}(b)=\frac{x}{\sqrt{x^2+b^2}}$. Consequently, $f'(x)>0$ and the function $f(x)$ is an increasing function for $a\leq b$. similarity, if $a>b$, then $f(x)$ is a decreasing function for every $x\geq 1$. 
\qed

\begin{lemma}\label{lemma3}
Let $g(x, y)=\sqrt{x^2+y^2}$ with $ y\geq 2$. Then $f(x, y)$ is an increasing function for every $x\geq 1$. 
\end{lemma}
\proof
We have $f'(x, y)=\frac{x}{\sqrt{x^2+y^2}}$. Since $x\geq 1$, $f'(x, y)>0$ and function $f(x, y)$ is an increasing function for every $x\geq 1$. 
\qed

\section{Extremal trees with the maximum Sombor index}
In this section, we characterize the extremal trees with maximum Sombor index among the trees with given degree sequence. We propose a technique to construct these trees. to do this, we first state some properties of a maximum optimal tree.\\

\begin{theorem}\label{theorem1}
Let $T$ be a maximum optimal tree with a path $v_0 v_1 v_2 \cdots v_k v_{k+1}$ in $T$, where $v_0$ and $v_{k+1}$ are leaves. For $i\leq \frac{t+1}{2}$ and $i+1\leq j \leq k-i+1$
\begin{enumerate}
\item[(i)] if $i$ is odd, then $d(v_i) \geq d(v_{k-i+1}) \geq d(v_j)$,
\item[(ii)] if $i$ is even, then $d(v_i) \leq d(v_{k-i+1}) \leq d(v_j)$.
\end{enumerate}
\end{theorem}

\proof
Let $T$ be a maximum optimal tree with the degree sequence $D$. We prove the result by induction on $i$. For $i=1$, we show that $d(v_1)\geq d(v_k)\geq d(v_j)$ where $2\leq j \leq k$. We suppose for contradiction that $d(v_1)<d(v_j)$ for some $2\leq j \leq k$. We consider a new tree $T'$ obtained from $T$ by changing edges $v_0 v_1$ and $v_j v_{j+1}$ to edges $v_0 v_j$ and $v_1 v_{j+1}$ such that no other edges are changed. Note that $T$ and $T'$ have the same degree sequence. Therefore, using Lemmas \ref{lemma1}-\ref{lemma3}, and since $d(v_{j+1})>1$, we have
\begin{align*}
SO(T')-SO(T)&=\sqrt{d(v_0)^2+d(v_j)^2}+\sqrt{d(v_1)^2+d(v_{j+1})^2}\\
&-\Big(\sqrt{d(v_0)^2+d(v_1)^2}-\sqrt{d(v_j)^2+d(v_{j+1})^2}\Big)\\
&=\Big(\sqrt{d(v_j)^2+1}-\sqrt{d(v_1)^2+1}\Big)\\
&+\Big(\sqrt{d(v_1)^2+d(v_{j+1})^2}-\sqrt{d(v_j)^2+d(v_{j+1})^2} \Big)\\
&=f(1)-f(d(v_{j+1}))>0,
\end{align*}

which is a contradiction with the maximum optimality $T$. Thus, $d(v_1)\geq d(v_j)$ for every $2\leq j \leq k$. similarity, we can get $d(v_1)\geq d(v_k)$ and $d(v_k)\geq d(v_j)$. Therefore, we have $d(v_1)\geq d(v_k)\geq d(v_j)$ where $2\leq j \leq k$. So, we suppose that the result holds for smaller values of $i$. \\
If $i\geq 2$ is even, then $i-1$ is odd and by the induction hypothesis, $d(v_{i-1}) \geq d(v_{k-i+1}) \geq d(v_j)$ for $i+1\leq j \leq k-i+1$. We suppose for contradiction that $d(v_i)>d(v_j)$ for some $i+1\leq j \leq k-i+1$. We consider a new tree $T''$ obtained from $T$ by changing edges $v_{i-1} v_i$ and $v_j v_{j+1}$ to edges $v_{i-1} v_j$ and $v_i v_{j+1}$ with the degree sequence $D$. Also, in tree $T''$, other edges are the same edges in tree $T$. \\
By the induction hypothesis, $d(v_{i-1})\geq d(v_{j+1})$. Therefore, by applying Lemma \ref{lemma2}, we have

\begin{align*}
SO(T'')-SO(T)&=\sqrt{d(v_{i-1})^2+d(v_j)^2}+\sqrt{d(v_i)^2+d(v_{j+1})^2}\\
&-\Big(\sqrt{d(v_{i-1})^2+d(v_i)^2}-\sqrt{d(v_j)^2+d(v_{j+1})^2}\Big)\\
&=\Big(\sqrt{d(v_j)^2+d(v_j)^2}-\sqrt{d(v_{i-1})^2+d(v_i)^2}\Big)\\
&+\Big(\sqrt{d(v_{j+1})^2+d(v_{i})^2}-\sqrt{d(v_{j+1})^2+d(v_{j})^2} \Big)\\
&=f(d(v_{i-1})-f(d(v_{j+1}))>0.
\end{align*}

 This contradiction with the maximum optimality of $T$. Therefore, $d(v_i)\leq d(v_j)$ for $i+1\leq j \leq k-i+1$. Similarity, we have $d(v_i)\leq d(v_{k-i+1})$ and $d(v_{k-i+1}) \leq d(v_j)$. Consequently, for $i$ even, $d(v_i) \leq d(v_{k-i+1}) \leq d(v_j)$ where $i+1\leq j \leq k-i+1$. For odd $i>2$, with similarity technique, we can get $d(v_i) \geq d(v_{k-i+1}) \geq d(v_j)$ for $i+1\leq j \leq k-i+1$. 
\qed

Suppose that $L_i$ denotes the set of vertices adjacent to the closet leaf at a distance $i$. Thus, $L_0$ and $L_1$ denote the set of leaves and the set of vertices that are adjacent to the leaves. Let $d^m=min\{d(u)~:~ u\in L_1 \}$ and $L_1^m$ be the set of leaves whose adjacent vertices have degree $d^m$ in $T$. We suppose that $\overline{L_1^m}$ denote the set of leaves $v$ such that $v \notin L_1^m$. \\
We construct a new tree $T'_i$ from tree $T$ and tree $T_i$ rooted at $v_i$ by identifying the root $v_i$ with a vertex $v \in L_1^m$. 

\begin{theorem}\label{theorem2}
Let $T'_1$ and $T'_2$ are obtained from $T$ by identifying the root $v_i$ of $T_i$ with $u'\in L_1^m$ and $v' \in \overline{L_1^m}$, respectively. Then, $SO(T'_1)\geq SO(T'_2)$. 
\end{theorem}
\proof
We suppose that $u$ and $v$ are adjacent to $u'$ and $v'$, respectively. Using Theorem \ref{theorem1}, $d(u)\leq d(v)$. Therefore, we have

\begin{align*}
SO(T'_1)-SO(T'_2)&=\sqrt{(d(v_{i})+1)^2+d(u)^2}+\sqrt{d(u)^2+1}\\
&-\Big(\sqrt{(d(v_{i})+1)^2+d(v)^2}-\sqrt{d(v)^2+1}\Big)\\
&=\Big(\sqrt{(d(v_{i})+1)^2+d(u)^2}-\sqrt{(d(v_{i})+1)^2+d(v)^2}\Big)\\
&+\Big(\sqrt{d(v)^2+1}-\sqrt{d(u)^2+1} \Big)\\
&=f(1)-f(d(v_{i})+1)>0.
\end{align*}

Therefore, $SO(T'_1)\geq SO(T'_2)$. 
\qed

We use a similar technique in \cite{Wang}, for constructing tree $T$ with a fixed degree sequence $D$ such that $T$ is the maximum optimal tree among the trees with degree sequence $D$. We propose the following algorithms to construct such trees. \\

\noindent {\bf Algorithm 1.} (Construction of subtrees)
\begin{enumerate}
\item Given the degree sequence of the non-leaf vertices as $D=\big(d_1, d_2, \ldots, d_m\big)$ in descending order. 
\item If $d_m \geq m-1$, then using Theorem \ref{theorem1}, the vertices with degrees $d_1, d_2, \ldots, d_{m-1}$ are in $L_1$. Tree $T$ produces by rotted at $u$ with $d_m$ children whose their degrees are $d_1, d_2, \ldots, d_{m-1}$ and $d_{m}-m+1$ leaves adjacent to $u$. 
\item If $d_m\leq m-2$, then we produce subtree $T_1$ by rotted at $u_1$ with $d_m-1$ children with degrees $d_1, d_2, \ldots, d_{d_{m-1}}$ such that $u_1 \in L_2$ and the children of $u_1$ are in $L_1$. \\
Subtree $T_2$ is constructed by rooted at $u_2$ with $d_{m-1}-1$ children whose degrees are $d_{d_m}, d_{d_m+1}, \ldots, d_{(d_m-1)+(d_{m-1}-1)}$. Then do the same to get subtrees $T_3, T_4, \ldots$ until $T_k$ satisfies the condition of step (2). In this case, we have $d(v_k)=d_{m-k+1}$.  
\end{enumerate}

\vspace{0.4 cm}

\noindent {\bf Algorithm 2.} (Merge of subtrees)
\begin{enumerate}
\item Set $T=T_i$ and $i=k$. We produce a new tree $T'_{i-1}$ from $T$ and $T_{i-1}$ rooted at $v_{i-1}$ by identifying the root $v_{i-1}$ with a vertex $v \in L_1^m$. Using Theorem \ref{theorem2}, tree $T'_{i-1}$ is a maximum optimal tree among trees with the same degree sequence. 
\item Consider $i=k-1, k-2, \ldots, 1$ and $T=T_i$. Tree $T'_{i-1}$ from $T$ and $T_{i-1}$ by the same method of step (1). We construct trees $T'_{k-2}, T'_{k-3}, \ldots, T'_1$. 
\item $T=T'_1$ is the maximum optimal tree with given degree sequence $D=\big(d_1, d_2, \ldots, d_m\big)$. 
\end{enumerate}

\begin{figure}[htp]
\begin{center}
\includegraphics[scale=0.75]{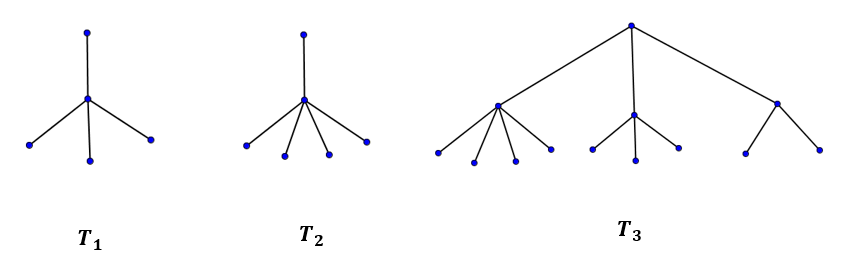}
\caption{Construction of subtrees using Algorithm 1}\label{f1}
\end{center}
\end{figure}

\vspace{0.25cm}

\begin{figure}[htp]
\begin{center}
\includegraphics[scale=0.75]{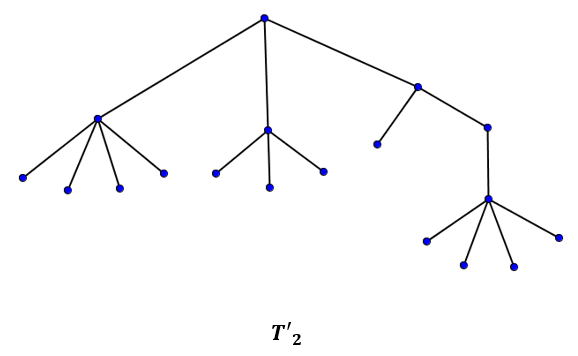}
\caption{Merge of subtrees using Algorithm 2}\label{f2}
\end{center}
\end{figure}

\vspace{0.25cm}

\begin{figure}[htp]
\begin{center}
\includegraphics[scale=0.75]{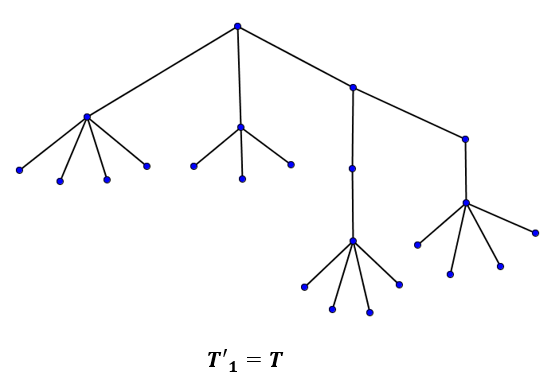}
\caption{A maximum optimal tree $T$ with degree sequence $\big( 5, 5, 5, 4, 3, 3, 2, 2\big)$.}\label{f3}
\end{center}
\end{figure}

\vspace{0.25cm}

\begin{example}{\em
In this example, we propose a maximum optimal tree with given degree sequence $D=\big( 5, 5, 5, 4, 3, 3, 2, 2\big)$. Using step (3) of Algorithm 1, we have subset $T_1$ with 1 child whose has degree 5. For new degree sequence $D_1=\big(5, 5, 4, 3, 3, 2\big)$, we construct tree $T_2$ and have new degree sequence $D=\big(5, 4, 3, 3\big)$ (Figure \ref{f1}). It is easily seen that $D_2$ satisfies the condition of step (2).\\
Using Algorithm 2, we attach subtrees $T_2$ to $T_3$ for constructing $T'_2$ (Figure \ref{f2}) and $T_2$ to $T'_2$ for constructing the maximum optimal tree $T'_1=T$ (Figure \ref{f3}). }
\end{example}

\vspace*{0.5 cm}
\noindent\textbf{Acknowledgements} The author would like to thank Professor Ivan Gutman for his useful comments and suggestions.

\end{document}